\documentclass[11pt]{article}
\usepackage{graphicx}
\usepackage{amsthm,amsmath,amsfonts,amssymb} 
\usepackage{color}
\catcode`\@=11 \@addtoreset{equation}{section}

\catcode`\@=12

\newtheorem{Theorem}{Theorem}[section]
\newtheorem{Proposition}{Proposition}[section]
\newtheorem{Lemma}{Lemma}[section]
\newtheorem{Corollary}{Corollary}[section]
\newtheorem{Remark}{Remark}[section]

\newcommand{\bTheorem}[1]{
\begin{Theorem} \label{T#1} }
\newcommand{\eT}{\end{Theorem}}

\newcommand{\bProposition}[1]{
\begin{Proposition} \label{P#1}}
\newcommand{\eP}{\end{Proposition}}

\newcommand{\bRemark}[1]{
\begin{Remark} \label{R#1}}
\newcommand{\eR}{\end{Remark}}

\newcommand{\bLemma}[1]{
\begin{Lemma} \label{L#1} }
\newcommand{\eL}{\end{Lemma}}

\newcommand{\bCorollary}[1]{
\begin{Corollary} \label{C#1} }
\newcommand{\eC}{\end{Corollary}}

\newcommand{\bFormula}[1]{
\begin{equation} \label{#1}}
\newcommand{\eF}{\end{equation}}

\newcommand{\sigmae}{\sigma_\lambda}
\newcommand{\we}{\vc{w}_\lambda}

\newcommand{\vr}{\varrho}
\newcommand{\vre}{\vr_\lambda}

\newcommand{\vue}{\vu_\lambda}

\newcommand{\vu}{\vc{u}}
\newcommand{\vc}[1]{{\bf #1}}
\newcommand{\Div}{{\rm div}_x}
\newcommand{\Grad}{\nabla_x}

\newcommand{\tn}[1]{\mbox {\F #1}}
\newcommand{\dx}{{\rm d} {x}}

\newcommand{\Ome}{\Omega_\lambda}

\newcommand{\intOe}[1]{\int_{\Omega_\lambda} #1 \ \dx}

\newcommand{\ep}{\varepsilon}

\font\F=msbm10 scaled 1000

\definecolor{Cgrey}{rgb}{0.85,0.85,0.85}
\definecolor{Cblue}{rgb}{0.50,0.85,0.85}
\definecolor{Cred}{rgb}{1,0,0}
\definecolor{fancy}{rgb}{0.10,0.85,0.10}

\newcommand\Cbox[2]{%
    \newbox\contentbox%
    \newbox\bkgdbox%
    \setbox\contentbox\hbox to \hsize{%
        \vtop{
            \kern\columnsep
            \hbox to \hsize{%
                \kern\columnsep%
                \advance\hsize by -2\columnsep%
                \setlength{\textwidth}{\hsize}%
                \vbox{
                    \parskip=\baselineskip
                    \parindent=0bp
                    #2
                }%
                \kern\columnsep%
            }%
            \kern\columnsep%
        }%
    }%
    \setbox\bkgdbox\vbox{
        \color{#1}
        \hrule width  \wd\contentbox %
               height \ht\contentbox %
               depth  \dp\contentbox
        \color{black}
    }%
    \wd\bkgdbox=0bp%
    \vbox{\hbox to \hsize{\box\bkgdbox\box\contentbox}}%
    \vskip\baselineskip%
}

\date{}

\providecommand{\com}{}

\begin{document}


\title{Robustness of strong solutions to the compressible Navier-Stokes system}

\author{Peter Bella \and Eduard Feireisl \thanks{Eduard Feireisl acknowledges the support of the GA\v CR (Czech Science Foundation) project P201-13-00522S in the framework of RVO: 67985840.} \and Bum Ja Jin \and Anton\' \i n Novotn\' y }


\maketitle

\bigskip
\centerline{Max Planck Institute for Mathematics in the Sciences}
\centerline{Inselstrasse 22, 04103 Leipzig, Germany}

\bigskip
\centerline{Institute of Mathematics of the Academy of Sciences of the Czech Republic}

\centerline{\v Zitn\' a 25, 115 67 Praha 1, Czech Republic}



\bigskip

\centerline{Department of Mathematics Education, Mokpo National University}

\centerline{Muan 534-729, South Korea}

\bigskip

\centerline{IMATH, EA 2134, Universit\' e du Sud Toulon-Var
BP 20132, 83957 La Garde, France}



\maketitle

\bigskip





\begin{abstract}
We consider the Navier-Stokes system describing the time evolution of a compressible barotropic fluid confined to a bounded spatial domain in the 3-D physical space, supplemented with the Navier's slip boundary conditions. It is shown that the class of global in time strong solutions
is robust with respect to small perturbations of the initial data. Explicit qualitative estimates are given also in terms of the
shape of the underlying physical domain, with applications to problems posed on thin cylinders.
\end{abstract}

\medskip

{\bf Key words:} Compressible Navier-Stokes system, regular solutions, robustness, thin cylindrical domain

\medskip


\section{Introduction}
\label{i}

In a series of very interesting papers \cite{ChCoRoTi}, \cite{MaroSa}, \cite{RoSad}, the authors discuss the hypothetical possibility of verifying regularity of
solutions to the incompressible Navier-Stokes system through a numerical test applied to a suitable finite set of initial data. The core of such a method is the property of \emph{robustness} of the class of strong solutions. More specifically, any small perturbation of the initial data
giving rise to a global-in-time strong solution enjoys the same property. In other words, the class of strong solutions is open in a suitable
topology. Another example of robustness is regularity of the solutions to the compressible Navier-Stokes system in the low Mach number regime
studied by Hagstrom and Lorenz \cite{HALO} in the situation that the limit (incompressible) system admits a strong solution. Last but not least,
robustness of the class of smooth solutions of the incompressible Navier-Stokes system plays a crucial role in the seminal paper
by Raugel and Sell \cite{RaSe2} concerning problems on thin domains.

Motivated by the previous results, we consider the \emph{compressible Navier-Stokes system}:

\Cbox{Cgrey}{

\bFormula{i1}
\partial_t \vr + \Div (\vr \vu) = 0,
\eF
\bFormula{i2}
\partial_t (\vr \vu) + \Div (\vr \vu \otimes \vu) + \Grad p(\vr) = \Div \tn{S}(\Grad \vu),
\eF

}

\noindent
where the unknowns are the mass density $\vr = \vr(t,x)$ and the velocity field $\vu = \vu(t,x)$, while
$p=p(\vr)$ is the pressure and $\tn{S}(\Grad \vu)$ the viscous stress obeying Newton's rheological law

\Cbox{Cgrey}{

\bFormula{i3}
\tn{S}(\Grad \vu) =
\mu \left( \Grad \vu + \Grad^t \vu - \frac{2}{3} \Div \vu
\tn{I} \right)  + \eta  \Div \vu \tn{I} ,
\eF

}

\noindent
with the shear viscosity coefficient
$\mu > 0$ and the bulk viscosity coefficient $\eta \geq 0$. Note that we have deliberately omitted the effect of external forces to simplify the presentation.

The fluid is confined to a bounded physical domain $\Omega \subset R^3$,
on the boundary of which we impose the slip conditions

\Cbox{Cgrey}{

\bFormula{i4}
\vu \cdot \vc{n}|_{\partial \Omega} = 0,\ \left[ \tn{S} (\Grad \vu) \cdot \vc{n} \right] \times \vc{n} |_{\partial \Omega} = 0,
\eF

}

\noindent
where $\vc{n}$ is the outer normal vector.
The motion originates from the initial state
\bFormula{i5}
\vr(0, \cdot) = \vr_{0}, \ \vu (0, \cdot) = \vu_{0}.
\eF

We suppose that the problem (\ref{i1} - \ref{i5}) admits a smooth solution $[\vr,\vu]$
defined on a time interval $[0,T]$. We consider a suitable perturbation $[\vr_{0,\lambda}, \vu_{0,\lambda}]$ of the initial data and suppose, in view of
future applications to thin and slim domains, that the domain $\Omega = \Omega_\lambda \subset R^3$ depends on the parameter
$\lambda = [\ep, d, V]$ in the following way:
\bFormula{i6-}
d = {\rm diam}[\Omega_\lambda], \ 0 < V = | \Omega_\lambda | < \infty;
\eF
there exists a uniformly $C^4-$domain $\mathcal{O}$ such that
\bFormula{i6}
\Ome = \ep \mathcal{O},\ 0 < \ep < 1.
\eF

Our aim is to identify a positive number $\omega = \omega(\lambda,T)$ such that the problem (\ref{i1} - \ref{i4}) admits a strong solution $[\vre, \vue]$ defined
on the time interval $[0,T]$ whenever the initial data $[\vr_{0,\lambda}, \vu_{0,\lambda}]$ satisfy
\bFormula{i7}
\left\| \vr_{0,\lambda} - \vr_0 \right\|_{W^{1,4}(\Ome)} + \left\| \vu_{0,\lambda} - \vu_0 \right\|_{W^{2,2}(\Ome;R^3)} < \omega(\lambda,T).
\eF
It is important that the specific form of $\omega$ depends solely on $T$, the geometric properties of the model domain $\mathcal{O}$, and
on the norm of certain derivatives of the smooth solution $[\vr, \vu]$. Applications will be given to problems posed on thin channels that motivate the choice of the slip condition (\ref{i4}), see \cite{BeFeNov}.

We consider the class of strong solutions compatible with (\ref{i7}) introduced by Cho, Choe and Kim \cite{ChoChoKi}. Note that a local existence result for the
problem (\ref{i1}-\ref{i4}) was also proved by Hoff \cite{Hoff2012} in a class of solutions enjoying better regularity than indicated by the
norms in (\ref{i7}). Adapting the results of the present paper to Hoff's class would require refined estimates accompanied with numerous technical difficulties.

The paper is organized as follows. In Section \ref{k}, we introduce a number of preliminary results concerning the scaled versions of
various embedding relations as well as the Lam{\' e} system associated with elliptic part of the momentum equation (\ref{i2}). In Section \ref{m}, we introduce the class of strong solutions to the problem (\ref{i1} - \ref{i4}) and state our main result concerning robustness. Then, in Section \ref{u}, we establish the necessary estimates on the perturbed solutions $[\vre, \vue]$ and prove the main result. Finally, Section \ref{A}
contains some applications of the main theorem
to problems on thin spatial cylinders.

\medskip

{\bf Notation:} By $C,C_1,\ldots$ we will denote positive constants which could depend on the parameters of the model domain $\mathcal{O}$, the smooth solution $[\vr, \vu]$, \emph{but neither} on $\lambda$ \emph{nor} on the initial data for the ``perturbed'' problem. The value of these constants could vary from expression to expression. By $a \lesssim b$ we mean that there exists constant $C$ such that $a \le C b$, similarly $a \backsimeq b$ means
$C_1 b < a < C_2b$.

\section{Preliminaries, useful inequalities, Lam{\' e} system}
\label{k}

In this section we will state (and prove) some auxiliary results, which will be handy later. To begin, we recall that
$\mathcal{O} \subset R^3$ is a uniformly $C^k-$domain, $k=0,1,\dots$ of type $(\alpha, \beta, K)$ if
for each point of $x_0 \in \partial \mathcal{O}$, there is a function $h \in C^{k}(R^2)$,
$\| h \|_{C^{k}(R^2)} \leq K$,
and
\[
U_{\alpha, \beta,h} = \{ (y, x_3) \ | h(y) - \beta < x_3 < h(y) + \beta ,\
|y| < \alpha \}
\]
such that, after suitable translation and rotation of the coordinate axes,
$x_0 = [0,0,h(0)]$ and
\[
\mathcal{O} \cap U_{\alpha, \beta, h} = \{ (y, x_3) \ | h(y) - \beta < x_3 < h(y) ,\
|y| < \alpha \},
\]
\[
\partial \mathcal{O} \cap U_{\alpha, \beta, h} =
\{ (y, x_3) \ |  x_3 = h(y) ,\
|y| < \alpha \},
\]
see Adams \cite{A}. All estimates listed below depend solely on the parameters $(\alpha, \beta,K)$ of the model domain $\mathcal{O}$.

\subsection{Korn's inequality}

{\

Korn's inequality on $\mathcal{O}$ reads
\[
\int_{\mathcal{O}} |\Grad \vc{v}|^2 \dx \le C(\alpha, \beta, K) \left( \int_{\mathcal{O}} |\vc{v}|^2\ \dx + \int_{\mathcal{O}} \tn{S} (\Grad \vc{v}) : \Grad \vc{v}\ \dx\right)
\]
for any $\vc{v} \in W^{1,2}(\mathcal{O}; R^3)$, see Dain \cite{Dain}, Reshetnyak \cite{Resh}.  Rescaling to $\Ome$ gives rise to
\bFormula{k1}
\ep^2 \intOe{ |\Grad \vc{v} |^2 } \leq C(\alpha, \beta, K) \left( \intOe{ |\vc{v}|^2 } + \ep^2 \intOe{ \tn{S} (\Grad \vc{v}) : \Grad \vc{v} }
\right)
\eF
for any $\vc{v} \in W^{1,2}(\Ome; R^3)$.

\bRemark{K1}

Note that the above estimates depend solely on the constants $\alpha$, $\beta$, and $K$ characterizing the shape of $\partial \mathcal{O}$ but not on
$|\mathcal{O}|$.

\eR

}

\subsection{Sobolev embedding}

Similarly to the preceding section, starting with Sobolev embedding relation
\[
\left\| v \right\|_{L^q(\mathcal{O})} \leq c(\alpha, \beta, K) \left( \left\| v \right\|_{L^p(\mathcal{O})} + \left\| \Grad v \right\|_{L^p(\mathcal{O};R^3)} \right) \ \mbox{for}\ 1 \leq p < 3, \ p \leq q \leq \frac{3p}{3 - p}
\]
we deduce
\bFormula{k2}
\left\| v \right\|_{L^q(\Ome)} \lesssim \ep^{3 \left( \frac{1}{q} - \frac{1}{p} \right)} \left( \left\| v \right\|_{L^p(\Ome)} + \ep \left\| \Grad v \right\|_{L^p(\Ome;R^3)} \right) \ \mbox{for}\ 1 \leq p < 3, \ p \leq q \leq \frac{3p}{3 - p}.
\eF
Note that the uniformly $C^1-$domains possess the $W^{1,p}-$extension property with the corresponding constant depending solely on the type
parameters $(\alpha, \beta, K)$.

Finally, by the same token,
\bFormula{k3}
\| v \|_{L^\infty(\Ome)} \lesssim \ep^{-\frac{3}{p}}\left( \left\| v \right\|_{L^p(\Ome)} + \ep \left\| \Grad v \right\|_{L^p(\Ome;R^3)} \right) ,\quad p>3,
\eF
for any $v \in W^{1,p}(\Ome)$.

\subsection{Gagliardo-Nirenberg inequality}

Combining (\ref{k2}) with the standard interpolation inequality
\[
\|v\|_{L^4}\lesssim \|v\|_{L^2}^{\frac{1}{4}}\|v\|_{L^6}^{\frac{3}{4}}
\]
we obtain
\bFormula{k5}
\Big(\intOe{ |v|^4 }\Big)^{\frac{1}{4}} \lesssim  \ep^{-\frac{3}{4}}\Big(\intOe{ |v|^2 }\Big)^{\frac{1}{8}}\left(\ep^2\intOe{ |\Grad v|^2}+\intOe{ |v|^2 }\right)^{\frac{3}{8}}
\eF
for any ${v} \in W^{1,2}(\Ome)$.

\subsection{Poincar\'e inequality}

Finally, using the same rescaling as above we get

{

\bFormula{k4}
\Big(\intOe{ |v-\bar{v}|^p }\Big)^{\frac{1}{p}} \lesssim d \left( \intOe{ |\Grad {v} |^p }\right)^{\frac{1}{p}},
\ \mbox{where}\ \bar{v} = \frac{1}{V} \intOe{ v }, \ \quad 1\leq p<\infty
\eF
for any ${v} \in W^{1,p}(\Ome)$.

}

\subsection{Lam{\' e} system}
\label{LM}

We consider the elliptic system
\bFormula{L1}
{\com \Div \tn{S}(\Grad \vc{w}) 
\equiv
\mu \Delta \vc{w} + \left( \frac{1}{3} \mu + \eta \right) \Grad \Div \vc{w}
= \vc{g} \ \mbox{in}\ \mathcal{O},}
\eF
supplemented with the boundary conditions
\bFormula{L2}
\vc{w} \cdot \vc{n}|_{\partial \mathcal{O}} = \left[ \left( \Grad \vc{w} + \Grad \vc{w}^t \right) \cdot \vc{n} \right] \times \vc{n}|_{\partial \mathcal{O}} = 0,
\eF
known as the \emph{Lam\' e system} in linear elasticity.

The standard elliptic theory developed by Agmon, Douglis, and Nirenberg \cite{ADN} (cf. also Hoff \cite[Lemma 2.2]{Hoff2012}) yields the {\it a priori} bound
\bFormula{L3}
\| \nabla^2_x \vc{w} \|_{L^p(\mathcal{O}; R^{3 \times 3})} \leq c(\alpha, \beta, K,p) \left( \| \vc{g} \|_{L^p(\mathcal{O};R^3)} +
\| \vc{w} \|_{L^p (\mathcal{O}; R^3)} \right), \ 1 < p < \infty,
\eF
which, after rescaling, gives rise to
\bFormula{L4}
\ep^2 \| \nabla^2_x \vc{w} \|_{L^p(\Ome; R^{3 \times 3})}
\lesssim \Big[ \ep^2\| \Div \tn{S}(\Grad \vc{w} ) \|_{\com L^p(\Ome; R^{3})} + \| \vc{w} \|_{L^p (\Ome; R^3)} \Big], \ 1 < p < \infty,
\eF
for any $\vc{w} \in W^{2,p}(\Ome;R^3)$ satisfying the boundary conditions(\ref{L2}) on $\partial \Ome$.

Finally, a similar argument implies
\bFormula{L5}
\ep \| \Grad \vc{w} \|_{L^p(\Ome; R^{3)}} \lesssim \Big[ \ep \| G \|_{L^p(\Ome)} + \| \vc{w} \|_{L^p (\Ome; R^3)} \Big], \ 1 < p < \infty,
\eF
for any $\vc{w} \in W^{1,p}(\Ome;R^3)$, $\vc{w} \cdot \vc{n}|_{\partial \Ome} = 0$, where
\[
\Grad G = \Div \tn{S} (\Grad \vc{w} ) \ \mbox{in the sense that}\
\intOe{ \tn{S} (\Grad \vc{w}) : \Grad \varphi } = \intOe{ G \Div \varphi }
\]
for all $\varphi \in W^{1,p}(\Ome;R^3)$, $\varphi \cdot \vc{n}|_{\partial \Ome} = 0$.

\section{Main results}
\label{m}

We start with the class of strong solutions introduced by Cho, Choe, and Kim \cite{ChoChoKi}, namely
\[
\vre \in C([0,T]; W^{1,q}(\Ome; R^3), \ \vue \in C([0,T]; W^{2,2}(\Ome; R^3)),
\]

The following result is an easy adaptation of \cite[Theorem 7]{ChoChoKi}:

\bProposition{M1}
Let $3 < q \leq 6$ be given and let the initial data $[\vr_{0,\lambda}, \vu_{0,\lambda}]$ belong to the class
\[
\vr_{0,\lambda} \in W^{1,q}(\Ome), \ \vu_{0,\lambda} \in W^{2,2}(\Ome; R^3),
\]
\[
\vr_{0,\lambda} \geq \underline{\vr} > 0 \ \mbox{in}\ \Ome, \ \vu_{0,\lambda} \cdot \vc{n}|_{\partial \Ome} =
[\tn{S} (\Grad \vu_{0,\lambda}) \cdot \vc{n}] \times \vc{n} |_{\partial \Ome} = 0.
\]
Let $p \in C^1[0, \infty)$ be an increasing function of the density.

Then there exists $0< T_{\rm max} \leq \infty$ such that the problem (\ref{i1} - \ref{i5}) admits a (unique) strong solution $[\vre, \vue]$
defined on a maximal time interval $[0, T_{\rm max})$,
\[
\vre \in C([0,T]; W^{1,q}(\Ome)), \ \vue \in C([0,T]; W^{2,2}(\Ome; R^3)),
\]
\[
\partial_t \vre \in  C([0,T]; L^2(\Ome)), \ \vue \in L^2(0,T; W^{2,q}(\Ome;R^3)), \ \partial_t \vue \in L^2 (0,T; W^{1,2}(\Ome;R^3))
\]
\[
\sqrt{\vre} \partial_t \vue \in L^\infty(0,T; L^2(\Ome;R^3))
\]
for any $0 < T < T_{\rm max}$.

If $T_{\rm max} < \infty$ then
\bFormula{M2}
\limsup_{t \to T_{\rm max}-} \left[ \| \vre(t, \cdot) \|_{W^{1,q}(\Ome)} + \| \vue (t, \cdot) \|_{W^{1,2}(\Ome;R^3)} \right] = \infty.
\eF

\eP

\bRemark{M1}

As a matter of fact, the local existence result \cite[Theorem 7]{ChoChoKi} was proved for the no-slip boundary conditions $\vu|_{\partial \Ome} = 0$. However,
in view of the regularity properties of solutions to the Lam\' e system stated in Section \ref{LM},
the method developed in \cite[Section 5]{ChoChoKi} can be easily adapted to the present setting.

\eR

\bRemark{M2}

Hoff \cite{Hoff2012} established a similar existence theory for the problem (\ref{i1} - \ref{i5}) for a class of data enjoying higher regularity, namely
$\vr_0 \in W^{2,2}(\Ome)$, $\vu_0 \in W^{3,2}(\Ome;R^3)$.

\eR

We are ready to state our main result.

\Cbox{Cgrey}{

\bTheorem{M1}
Let $\mathcal{O} \subset R^3$ be a uniformly $C^4-$bounded domain, $\Ome = \ep \mathcal{O}$, $0 <\ep \leq 1$,
\bFormula{MD1}
{\rm diam} [\Ome] = d < d_0, \ 0 < | \Ome | = V \leq V_0 < \infty.
\eF
Let
\bFormula{MP1}
p \in C^1[0, \infty) \cap C^2(0,\infty),\ p(0) = 0, \ p' (\vr) > 0 \ \mbox{for all}\ \vr > 0.
\eF
Suppose that the problem (\ref{i1} - \ref{i4})
admits a smooth (classical) solution $[\vr, \vu]$ in $[0,T] \times \Ome$, emanating from the initial data $[\vr_0, \vu_0]$.

Then there exists a constant $C$, depending only on $d_0$, $V_0$, the type parameters $(\alpha, \beta, K)$ of the model domain $\mathcal{O}$, and on
\[
U = \max \left\{ \com \| \vr \|_{L^\infty (\Ome\times (0,T))} , \| \vr^{-1} \|_{L^\infty(\Ome\times (0,T))} , \| \partial^j_t \vr \|_{W^{1, \infty}(\Ome\times (0,T))} , \| \partial^j_t \vu \|_{W^{2, \infty}(\Ome\times (0,T); R^3)},\ j = 1,2  \right\}
\]
such that for any initial data $[\vr_{0,\lambda}, \vu_{0,\lambda}]$,
\[
\vr_{0,\lambda} \in W^{1,4}(\Ome), \ \vu_{0,\lambda} \in W^{2,2}(\Ome; R^3),
\]
\[
\vr_{0,\lambda} \geq \underline{\vr} > 0 \ \mbox{in}\ \Ome, \ \vu_{0,\lambda} \cdot \vc{n}|_{\partial \Ome} =
[\tn{S} (\Grad \vu_{0,\lambda}) \cdot \vc{n}] \times \vc{n} |_{\partial \Ome} = 0,
\]
satisfying
\bFormula{MM1}
\left\| \vr_{0,\lambda} - \vr_0 \right\|_{W^{1,4}(\Ome)} + \left\| \vu_{0,\lambda} - \vu_0 \right\|_{W^{2,2}(\Ome;R^3)} \leq \omega(\lambda,T),\
\lambda \equiv [d, V, \ep],
\eF
\bFormula{MM2}
\omega(\lambda, T) = \exp \left[ - C \left( {\ep^{-16/5}} +  V^{-1/4} \ep^{-1}
 \right) T \right] \min \left\{ \com \ep^{5}, \ep^{3/2}V^{1/2} \right\}.
\eF
the problem (\ref{i1} - \ref{i4}) admits a (unique) strong solution $[\vre, \vue]$ (in the sense specified in Proposition \ref{PM1}) in $[0,T] \times \Ome$,
\[
\vre(0, \cdot) = \vr_{0,\lambda}, \ \vue(0, \cdot) = \vu_{0, \lambda}.
\]
\eT

}

\bRemark{CLASS}

It is important to notice that we deal with perturbations of a \emph{general} smooth solution that
may not be stable. There are classical results by Matsumura and Nishida \cite{MANI}, Matsumura and
Padula \cite{MAPAD}, or, more recently, Valli and Zajaczkowski \cite{VAZA}, where the authors study
perturbations of a \emph{stable equilibrium solution}.

\eR 

\bRemark{MMM1-}

{\com It is important to note that the condition $\Omega = \ep \mathcal{O}$ represents only a restriction on the shape of the boundary of 
$\partial \Omega_\lambda$ 
expressed in terms of the parameters $[\alpha, \beta, K]$ but is independent of the volume and other features of the model domain $\mathcal{O}$. }

\eR

\bRemark{MMM1}

It is essential in the proof of Theorem \ref{TM1} that the
classical solution $[\vr, \vu]$ enjoys better smoothness properties than the solutions $[\vre, \vue]$. In particular, all derivatives appearing in the equations are assumed to be continuous and bounded up to the boundary, cf. the estimates in Section \ref{EI}.

\eR

The next section is devoted to the proof of Theorem \ref{TM1}. Applications to problems on thin domains are discussed in Section \ref{A}.

\section{Uniform estimates, proof of Theorem \ref{TM1}}
\label{u}

In view of Proposition \ref{PM1}, the local solution $[\vre, \vue]$ originating from the data $[\vr_{0,\lambda}, \vu_{0,\lambda}]$ can be continued up to the desired
time $T$ as long as we control the norms specified in (\ref{M2}). To this end, we derive a series of estimates valid on the existence interval $[0, T_{\rm max})$ of the solution $[\vre, \vue]$.

Introducing the differences $\sigmae = \vre - \vr$, $\we = \vue - \vc{u}$, we observe that
\bFormula{u1}
\partial_t \sigmae + \Div (\sigmae \vue) = - \Div ( \vr \we ),
\eF
with the associated renormalized equations (cf. DiPerna and Lions \cite{DL}),
\begin{eqnarray}
\label{u2}
\nonumber
 \partial_t H(\vre) + \Div \left( H(\vre) \vue \right) &=& - p(\vre) \Div \vue  \\
 \partial_t H(\vr) + \Div \left( H(\vr) \vc{u} \right) &=& - p(\vr) \Div \vc{u} ,
\end{eqnarray}
where the potential $H(\vr)$ is defined (modulo a linear function) through
\[
H'(\vr) \vr - H(\vr) = p(\vr).
\]

Similarly, subtracting the momentum equations we get
\bFormula{u3}
\vre \left( \partial_t \we + \vue \cdot \Grad \we \right) - \Div \tn{S} (\Grad \we) =
\Grad \left( p(\vr) - p(\vre) \right) - \sigmae \partial_t \vc{u} - \left( \vre \vue - \vr \vc{u} \right) \cdot \Grad \vc{u},
\eF
with the slip boundary conditions
\bFormula{u4}
\we \cdot \vc{n}|_{\partial \Ome} = 0, \ \left[ \tn{S} (\Grad \we) \cdot \vc{n} \right] \times \vc{n}|_{\partial \Ome} = 0,
\eF
and the initial conditions
\bFormula{u5}
\sigmae (0, \cdot) = \vr_{0,\lambda} - \vr_0, \ \we (0, \cdot) = \vu_{0,\lambda} - \vc{u}_0.
\eF

\subsection{Energy estimates}

Taking the scalar product of (\ref{u3}) with $\we$ and integrating the resulting expression over $\Ome$ we obtain
\bFormula{u6}
\frac{{\rm d}}{{\rm d}t} \intOe{ \frac{1}{2} \vre |\we|^2 } + \intOe{ \tn{S}\left(\Grad \we \right) : \Grad \we } \eF
\[
= \intOe{\left[ (p(\vre) - p(\vr) \Div \we  - \left( \sigmae \partial_t \vc{u} + \left( \vre \vue - \vr \vc{u} \right) \cdot \Grad \vc{u}
\right) \cdot \we \right] }.
\]
Next, the renormalized equations (\ref{u2}) give rise to
\bFormula{u7}
\frac{{\rm d}}{{\rm d}t} \intOe{ \left[ H(\vre) - H(\vr) \right] } = \intOe{ \left( p(\vr) \Div \vc{u} - p(\vre) \Div \vue \right) }.
\eF

Finally, since $\vr$ is strictly positive, we are allowed to multiply (\ref{u1}) by $H'(\vr)$ to deduce that
\bFormula{u8}
\frac{{\rm d}}{{\rm d}t} \intOe{ H'(\vr) \sigmae } = \intOe{ \left(
\partial_t H'(\vr) \sigmae + H'(\vr) \partial_t \sigmae \right) }
\eF
\[
= \intOe{ \left( \vre \Grad H'(\vr) \cdot \we - \sigmae p'(\vr) \Div \vc{u} \right) }.
\]

Summing up (\ref{u6} - \ref{u8}) we may infer that
\bFormula{u9}
\frac{{\rm d}}{{\rm d}t} \intOe{ \left[ \frac{1}{2} \vre |\we|^2  +
H(\vre) - H'(\vr) \sigma_\ep - H(\vr) \right] } + \intOe{ \tn{S}\left(\Grad \we \right) : \Grad \we }
\eF
\[
=-\intOe{ \left[ \we \cdot\Big(\sigmae \partial_t \vc{u}+\Grad \vc{u} \cdot (\vre \vue -\vr \vc{u} )+\sigmae \Grad H'(\vr) \right)
+\left(p(\vre)-p(\vr)-p'(\vr)\sigmae \Big)\Div \vc{u}\right] }.
\]

\bRemark{U1}

Formula \ref{u9} coincides with the so-called relative entropy inequality derived in \cite{FeNoJi}.

\eR

\subsection{Higher order estimates}

To deduce higher order estimates, we start by taking the scalar product of the momentum equation (\ref{u3}) with $\partial_t \we$ to obtain
\bFormula{u10}
\frac{{\rm d}}{{\rm d}t} \frac{1}{2} \intOe{ \tn{S} (\Grad \we ) : \Grad \we } + \intOe{ \vre |\partial_t \we |^2 }
\eF
\[
= -\intOe{ \left(
\vre  \Grad \we \cdot \vue + \vc{K}_\lambda \right) \cdot \partial_t \we  }
+\intOe{ (p(\vre)-p(\vr))\Div \partial_t \we  },
\]
where we have set
\[
\vc{K}_\lambda = \sigmae \partial_t \vc{u} + \Grad \vc{u} \cdot \left( \vre \vue - \vr \vc{u} \right).
\]

Now we differentiate equation (\ref{u3}) with respect to $t$:
\[
\vre \Big( \partial_t [\partial_t \we] + \vue \cdot \Grad [\partial_t \we ] \Big) - \Div \tn{S}(\partial_t \Grad \we )
\]
\[
= \Grad \partial_t \left( p(\vr) - p(\vre) \right) - \partial_t \vc{K}_\lambda - \partial_t \vre \Big( \partial_t \we + \vue \cdot \Grad \we \Big)
- \vre \Grad \we \cdot \partial_t \vue.
\]
Multiplying the above equation by $\partial_t \we$ and integrating over $\Ome$ we get
\bFormula{u11}
\frac{{\rm d}}{{\rm d}t} \frac{1}{2}
\intOe{ \vre
| \partial_t \we |^2 } + \intOe{ \tn{S} (\partial_t \Grad \we ) : \partial_t \Grad \we }
\eF
\[
=
\intOe{ \left( \partial_t p(\vre )- \partial_t p(\vr) \right) \Div \partial_t \we }
-\intOe{  \Big( \partial_t \vc{K}_\lambda + \partial_t \vre(\partial_t \we +( \vue \cdot \Grad )\we )+\vre \partial_t \vue \cdot \Grad
 \we  \Big) \cdot \partial_t \we  }.
\]

Finally, we multiply \eqref{u10} by $\ep^2$ and \eqref{u11} by $\ep^4$,
and, adding the resulting expression to \eqref{u9}, we derive
\bFormula{u12}
\frac{{\rm d}}{{\rm d}t} \intOe{ \left[ \frac{1}{2} \left( \vre |\we|^2 + \ep^4 \vre |\partial_t \we |^2 +\ep^2 \tn{S} (\Grad \we) : \Grad \we \right) + H(\vre) - H'(\vr) \sigmae - H(\vr)
\right] }
\eF
\[
+ \intOe{ \left[ \ep^2 \vre |\partial_t \we |^2 + \tn{S}(\Grad \we): \Grad \we +\ep^4 \tn{S}(\partial_t \Grad \we ) :
\partial_t \Grad \we  \right]  } =\sum_{j=1}^7 I_{j,\lambda},
\]
where
\begin{align*}
 I_{1,\lambda} &=  -\intOe{ \left[ \we \cdot \Big( \sigmae \partial_t \vc{u} +(\vre \vue - \vr \vc{u} )\cdot \Grad \vc{u} +\sigmae \Grad H'(\vr) \Big)
+\Big( p(\vre)-p(\vr)-p'(\vr)\sigmae \Big)\Div \vc{u} \right] }
,
\\
I_{2,\lambda} &= -\ep^2\intOe{
\vre ( \Grad \we \cdot \vue )\cdot \partial_t \we },
\\ I_{3,\lambda} &= - \ep^2\intOe{ \vc{K}_\lambda \cdot \partial_t \we },
\\ I_{4,\lambda} &= \ep^2\intOe{ \Big( p(\vre )-p(\vr) \Big) \Div \partial_t \we },
\\ I_{5,\lambda} &= \ep^4\intOe{ \Big( \partial_t p(\vre)-\partial_t p(\vr))\Div \partial_t \we },
\\ I_{6,\lambda} &= -\ep^4\intOe{  \left[ \partial_t \vre \Big( \partial_t \we +( \vue \cdot \Grad ) \we \Big)+\vre \partial_t \vue \cdot \Grad \we \right] \cdot \partial_t \we },
\\ I_{7,\lambda} &= -\ep^4\intOe{  \partial_t \vc{K}_\lambda \cdot \partial_t \we }.
\end{align*}

\subsection{Closing the estimates}

Denoting
\[
E_\lambda(t) = \intOe{ \left[ \frac{1}{2} \left( \vre |\we |^2 + \ep^4\vre |\partial_t \we |^2 +\ep^2 \tn{S} (\Grad \we ) : \Grad \we \right) + H(\vre) - H'(\vr) \sigmae - H(\vr) \right]{\com (t,x)} }
\]
and
\[
D_\lambda(t) =
\intOe{ \left[\ep^2 \vre |\partial_t \we |^2 + \tn{S}(\Grad \we): \Grad \we +\ep^4 \tn{S}(\partial_t \Grad \we ) :
\partial_t \Grad \we   \right] {\com (t,x)} },
\]
we aim to use (\ref{u12}) to close the estimates via a Gronwall type argument.

{ Since the density $\vr$ is bounded below and above uniformly for $t \in [0,T]$, we get
\[
\com C_1 - \|\sigmae\|_{L^\infty(\Ome)\times[0,T]} \leq \vre \le C_2 + \|\sigmae\|_{L^\infty(\Ome)\times[0,T]},
\]
and, consequently,
\bFormula{u13-1}
\left(C_1-\| \sigmae \|_{L^\infty(\com \Ome\times[0,T])}\right)\int_{\Ome}| \we|^2 \dx\leq  \int_{\Ome}\vre | \we |^2 \dx\leq \left(C_2+\| \sigmae \|_{L^\infty(\com\Ome\times[0,T])}\right)\int_{\Ome}| \we |^2 \dx,\eF
and
\bFormula{u13-2}
\left(C_1- \| \sigmae \|_{L^\infty(\com\Ome\times[0,T])}\right)\int_{\Ome}| \partial_t\we |^2 \dx\leq  \int_{\Ome}\vre |\partial_t \we |^2 \dx\leq \left(C_2+\| \sigmae \|_{L^\infty(\com\Ome\times[0,T])}\right)\int_{\Ome}|\partial_t \we |^2 \dx.
\eF

Moreover, Korn's inequality (\ref{k1}) implies
 \[
 \int_{\Ome}\ep^2|\Grad\we |^2 \dx \lesssim \int_{\Ome}\ep^2 \tn{S}(\Grad \we ) \ : \Grad \we  + |\we |^2 \dx
 \]
 and
\[
\int_{\Ome} \ep^4|\Grad \partial_t \we |^2\dx \lesssim \int_{\Omega_\ep}\ep^4 \tn{S}(\partial_t \Grad \we ) \ : \partial_t \Grad \we  +\ep^2 |\partial_t \we |^2 \dx.
\]

To simplify further discussion we shall assume that
\begin{equation}\label{h1}
 \| \sigmae (t,\cdot) \|_{L^\infty(\Ome)} \le { \frac{C_1}{2} },
\
\| \we (t,\cdot) \|_{L^\infty(\Ome)} \le 1
\end{equation}
for $0 \leq t \leq \min \{ T_{\rm max},T \}$. Such an inequality can be certainly achieved for $t = 0$ by a suitable choice of the initial data; whence it can be extended to the whole interval $[0, T_{\rm max}]$,
diminishing $\com T_{{\rm max}}$ as the case may be. We will later justify these assumptions, meaning we show that they hold on the whole interval
$[0,T]$ as soon as $[\vr_{0,\lambda}, \vu_{0,\lambda}]$ is sufficiently close to $[\vr_0, \vu_0]$.
Consequently, in accordance with the structural properties of $p$ stated in (\ref{MP1}), we may infer that
\bFormula{u13}
\intOe{ |\sigmae |^2 } \backsimeq \intOe{ H(\vre) - H'(\vr) \sigmae - H(\vr) }.
\eF

Combining the previous estimates we obtain
}
\begin{equation}\nonumber
\intOe{ \left[ |\we |^2 + \ep^2 |\Grad \we |^2 + |\sigmae |^2 +  \ep^4 |\partial_t \we |^2 \right](t, \cdot) } \lesssim E_\lambda (t)
\end{equation}
and
\begin{equation}\nonumber
\intOe{ \ep^2 |\partial_t \we |^2 + \ep^4 | \partial_t \Grad \we |^2(t, \cdot) } \lesssim D_\lambda(t);
\end{equation}
whence
\begin{align}
 \| \we (t,\cdot) \|_{L^2(\Ome)} &\lesssim \sqrt{E_\lambda (t)}, \label{est1}
\\ \| \Grad \we (t,\cdot) \|_{L^2(\Ome)} &\lesssim \ep^{-1} \sqrt{E_\lambda (t)}, \label{est2}
\\
 \| \sigmae (t,\cdot) \|_{L^2(\Ome)} &\lesssim \sqrt{E_\lambda (t)}\label{est3},
\\
\| \partial_t \we (t,\cdot) \|_{L^2(\Ome)} &\lesssim \ep^{-2} \sqrt{E_\lambda (t)}, \label{est4}
\\ \| \partial_t \we (t,\cdot) \|_{L^2(\Ome)} &\lesssim \ep^{-1} \sqrt{D_\lambda (t)}, \label{est5}
\\ \| \partial_t \Grad \we (t,\cdot) \|_{L^2(\Ome)} &\lesssim \ep^{-2} \sqrt{D_\lambda (t)}. \label{est6}
\end{align}

Furthermore, by virtue of Sobolev inequality (\ref{k2}),  we obtain
\begin{align}
\|\we (t,\cdot)\|_{L^p(\Ome)} &\lesssim \ep^{\frac{3}{p}-\frac{3}{2}}\sqrt{E_\lambda (t)}\label{est7},
\\ \|\partial_t\we (t,\cdot)\|_{L^p(\Ome)} &\lesssim \ep^{\frac{3}{p}-\frac{5}{2}}\sqrt{D_\lambda(t)} \label{est8}
\end{align}
for any $2 \le p \le 6$.

{
Next, using the Gagliardo-Nirenberg inequality (\ref{k5}) with \eqref{est4}, \eqref{est6}, we get
\begin{equation}
\begin{aligned}\label{est9}
\|\partial_t\we (t,\cdot)\|_{L^4(\Ome)} &\lesssim \ep^{-3/4}\|\partial_t \we (t,\cdot)\|_{L^2(\Ome)}^{1/4} \left( \ep \| \Grad \partial_t \we (t,\cdot)\|_{L^2(\Ome)} + \|\partial_t \we (t,\cdot)\|_{L^2(\Ome)}\right)^{3/4}
\\ &{\lesssim} \ep^{-3/4} \left( \ep^{-2} \sqrt{E_\lambda(t)} \right)^{1/4} \left( \ep^{-1} \sqrt{D_\lambda(t)} \right)^{3/4} = \ep^{-2} E_\lambda (t)^{\frac{1}{8}}D_\lambda (t)^{\frac{3}{8}}.
\end{aligned}
\end{equation}

Finally, application of Poincar\'e inequality (\ref{k4}) and H\"older's inequality, together with hypothesis (\ref{MD1}), gives rise to
\begin{equation}\label{est10}
\begin{aligned}
\|\sigmae (t,\cdot)\|_{L^4(\Ome)} &\lesssim d \|\Grad \sigmae (t,\cdot)\|_{L^4(\Ome)}+\|\overline{\sigmae (t,\cdot)}\|_{L^4(\Ome)}
\\ &\lesssim  d \|\Grad \sigmae (t,\cdot)\|_{L^4(\Ome)}+|\Ome|^{-3/4} \|\sigmae (t,\cdot)\|_{L^1(\Ome)}
\leq d \|\Grad \sigmae (t,\cdot)\|_{L^4(\Ome)}+|\Ome|^{-1/4} \|\sigmae (t,\cdot)\|_{L^2(\Ome)}
\\
&\lesssim d \|\Grad \sigmae (t,\cdot)\|_{L^4(\Ome)}+ V^{-1/4}    \sqrt{E_\lambda(t)}.
\end{aligned}
\end{equation}

}

\subsubsection{Estimating the integrals}
\label{EI}

In this section we will estimate the right hand side of \eqref{u12} (integrals $I_{1,\ep}$ through $I_{7,\ep}$). To perform the estimates we will often use that $[\rho,\vc{u}]$ is a smooth solution and as such all derivatives of $\rho$ and $\vc{v}$ are bounded in $L^\infty$, both in time and space. Moreover, we will also use that $\rho$ is bounded below away from zero in $[0,T]$.

\begin{enumerate}
\item
Applying H\"older's inequality several times and using \eqref{h1} we get
\begin{multline*}\nonumber
\left| \intOe{ \we \cdot \left( \sigmae \partial_t \vc{u} +(\vr \we + \sigmae \vc{u} + \sigmae \we )\cdot \Grad \vc{u} +\sigmae \Grad H'(\vr) \right)} \right|
\\ \lesssim \|\we \|_{L^2(\Ome;R^3)} \|\sigmae \|_{L^2(\Ome)} + \|\we \|_{L^2(\Ome;R^3)}^2 + \|\sigmae \|_{L^2(\Ome)}^2 
{\lesssim} E_\lambda 
\end{multline*}
and
\begin{equation*}
 \left| \int_{\Ome} \left( p(\vre)-p(\vr)-p'(\vr)\sigmae \right)\Div \vc{u} \right| \lesssim \|\sigmae \|_{L^2(\Ome)}^2 \lesssim E_\lambda,
\end{equation*}
where we have exploited \eqref{est1} and \eqref{est3}.

Adding these two inequalities together implies
\begin{equation}\label{estI1}
|I_{1,\lambda }| \lesssim E_\lambda .
\end{equation}

\item 
By H\"older's inequality, and the estimates \eqref{h1}, \eqref{est2}, \eqref{est7}, \eqref{est8}, we get
{ \begin{equation}\label{estI2}
 \begin{aligned}
 |I_{2,\lambda}| &= \ep^2 \left| \intOe{
\vre ( \Grad \we \cdot \vue )\cdot \partial_t \we } \right|
\\ &\le \ep^2 \| \vr + \sigmae \|_{L^\infty(\Ome)}
 \|\Grad \we \|_{L^2(\Ome;R^3)} \| \vu + \we \|_{L^\infty(\Ome;R^{3 \times 3})} \| \partial_t \we
\|_{L^2(\Ome;R^3)}
\\ &\lesssim \ep^2 \left(\ep^{-1} \sqrt{E_\lambda} \right) 
 \left( \ep^{-1} \sqrt{D_\lambda} \right) =
D^{1/2}_{\lambda} E^{1/2}_{\lambda}.
 \end{aligned}
\end{equation}}

\item
Since $\vc{K}_\lambda = \sigmae \partial_t \vc{u} + \Grad \vc{u} \cdot \left( \vr \we + \sigmae \vc{u} + \sigmae \we \right)$,
we may use \eqref{h1} to obtain
\[
\|\vc{K}_\lambda \|_{L^2(\Ome;R^3)} \lesssim \|\sigmae \|_{L^2(\Ome)} + \|\we \|_{L^2(\Ome;R^3)};
\]
whence, by virtue of \eqref{est1}, \eqref{est3}, and \eqref{est4}, we have
\begin{equation}\label{estI3}
|I_{3,\lambda}| \le \ep^2\|K_\lambda \|_{L^2(\Ome)}\|\partial_t\we \|_{L^2(\Ome;R^3)}
\lesssim \ep^2 (\|\sigmae \|_{L^2(\Ome)}+\|\we \|_{L^2(\Ome;R^3)})\|\partial_t\we \|_{L^2(\Ome;R^3)}
\lesssim E_\lambda.
\end{equation}
Similarly, using \eqref{est3} and \eqref{est6} we get
\begin{equation}\label{estI4}
|I_{4,\lambda}| \leq \ep^2\|p(\vre)-p(\vr)\|_{L^2(\Ome) }\|\Div \partial_t\we \|_{L^2(\Ome)}\lesssim
 \ep^2\|\sigmae \|_{L^2(\Ome)}\|\partial_t\Div \we \|_{L^2(\Ome) }\lesssim {E_\lambda^{1/2}}{D_\lambda^{1/2} }.
\end{equation}

\item
To estimate $I_{5,\lambda}$ we first observe that continuity equation \eqref{i1} with (\ref{u2}) imply that
\[
\partial_t(p(\vre)-p(\vr))= (\Grad p(\vr) - \Grad p(\vre ))\vc{u} - \Grad p(\vre ) \we + (p'(\vr)\vr - p'(\vre)\vre) \Div \vc{u} -
p'(\vre )\vre \Div \we,
\]
therefore, by \eqref{h1},
\begin{align*}
\|\partial_t(p(\vre)-p(\vr))\|_{L^2(\Ome)}
&\lesssim \| \Grad(p(\vre)-p(\vr))\|_{L^2(\Ome)}+ \|\we \|_{L^2(\Ome;R^3)}+\|\sigmae \|_{L^2(\Ome)}+\|\Grad \we \|_{L^2(\Ome; R^{3 \times 3})}
\\ &\lesssim \| \Grad \sigmae \|_{L^2(\Ome;R^3)}+\|\sigmae \|_{L^2(\Ome)}+ \|\we \|_{L^2(\Ome)}+\|\Grad \we \|_{L^2(\Ome;R^{3 \times 3})}.
\end{align*}
Since
\[
\|\Grad \sigmae \|_{L^2(\Ome;R^3)}\leq |\Ome|^{\frac{1}{4}}\| \Grad \sigmae \|_{L^4(\Ome;R^3) } =  V^{\frac{1}{4}} \| \Grad \sigmae \|_{L^4(\Ome;R^3)},
\]
the previous estimate together with (\ref{est1} - \ref{est3}), and \eqref{est6} imply
\bFormula{estI5}
|I_{5,\lambda}| \le \ep^4\|\partial_t(p(\vre)-p(\vr))\|_{L^2}\|\Div \partial_t\we \|_{L^2(\Ome)}
\lesssim \ep {E^{1/2}_\lambda}{D^{1/2}_\lambda }+\ep^{2} V^{\frac{1}{4}} \|\Grad \sigmae \|_{L^4(\Ome;R^3)}{D_\lambda^{1/2}}.
\eF

\item

 As
 \[
 { \partial_t\vre= -\vre \mbox{div}(\we+\vc{u})-\vue \cdot \Grad (\sigmae+\vr)},
 \]
  we may use H\"older's inequality together with \eqref{h1} to obtain
\[
\|\partial_t\vre\|_{L^2 (\Ome)}\lesssim |\Ome|^{\frac{1}{2}} + \|\Grad \we \|_{L^2(\Ome;R^{3 \times 3})}+\|\Grad \sigmae \|_{L^2(\Ome;R^3)}
\]
\[
\leq V^{ \com \frac{1}{4} } \left( V^{\com \frac{1}{4}} + \|\Grad  \we \|_{L^4(\Ome;R^{ 3 \times 3}) } + \|\Grad \sigmae \|_{L^4(\Ome; R^3)} \right).
\]
{Here and hereafter, we observe that
\[
\|\vue \|_{L^p(\Ome; R^3)} = \|\we + \vc{u} \|_{L^p(\Ome; R^3 )} \lesssim \|\we \|_{L^p(\Ome;R^3 )}+ V^{\frac{1}{p}}
\]
and, similarly,
\[
\|\vre\|_{L^p(\Ome)} = \|\vr + \sigmae \|_{L^p(\Ome )} \lesssim \|\sigmae \|_{L^p(\Ome)}+ V^{\frac{1}{p}}
\]
for any $1 \leq p<\infty$.
}

Thus the previous estimate together with \eqref{h1}, \eqref{est2}, and \eqref{est9} imply
\bFormula{estI6}
| I_{6, \lambda} | \leq \ep^4 \|\partial_t\we \|_{L^4(\Ome; R^3)} \|\partial_t\vre\|_{L^2 (\Ome)} \Big[
\|\partial_t\we \|_{L^4 (\Ome; R^3)} + \|\vue \|_{L^\infty(\Ome;R^3)}\|\Grad \we \|_{L^4(\Ome; R^{3 \times 3})} \Big]
\eF
\[
+ {\com \ep^4} \|\partial_t\we \|_{L^4(\Ome; R^3)}\|\vre \|_{L^\infty(\Ome)}\|\partial_t\vue \|_{L^4(\Ome;R^3)}\|\Grad \we \|_{L^2(\Ome; R^{3 \times 3})}
\]
\[
\lesssim \ep^4 \left( \ep^{-2} D_\lambda^{3/8} E_\lambda^{1/8} \right) V^{ \frac{1}{\com 4} } \left( V^{\frac{1}{\com 4}} + \|\Grad  \we \|_{L^4(\Ome;R^{ 3 \times 3}) } + \|\Grad \sigmae \|_{L^4(\Ome; R^3)} \right) \left( \ep^{-2} D_\lambda^{3/8} E_\lambda^{1/8} \right)
\]
\[
+ \ep^4 \left( \ep^{-2} D_\lambda^{3/8} E_\lambda^{1/8} \right) V^{ \frac{1}{\com 4} } \left( V^{\frac{1}{\com 4}} + \|\Grad  \we \|_{L^4(\Ome;R^{ 3 \times 3} } + \|\Grad \sigmae \|_{L^4(\Ome; R^3)} \right)  \|\Grad \we \|_{L^4(\Ome; R^{3 \times 3})}
\]
\[
+ \ep^4\left( \ep^{-2} D_\lambda^{3/8} E_\lambda^{1/8} \right) \left( \ep^{-2} D_\lambda^{3/8} E_\lambda^{1/8} + V^{1/4} \right) \left( \ep^{-1} \sqrt{E_\lambda} \right)
\]
\[
\lesssim V^{\frac{1}{\com 4}} D_\lambda^{3/4} E_\lambda^{1/4}  \left( {\com V^{\frac{1}{4}}} + \|\Grad  \we \|_{L^4(\Ome;R^{ 3 \times 3}) } + \|\Grad \sigmae \|_{L^4(\Ome; R^3)} \right)
\]
\[
+
\ep^2 V^{\frac{1}{\com 4}} D_\lambda^{3/8} E_\lambda^{1/8}  \left( {\com V^{\frac{1}{4}}} + \|\Grad  \we \|_{L^4(\Ome;R^{ 3 \times 3}) } + \|\Grad \sigmae \|_{L^4(\Ome; R^3)} \right)
\|\Grad \we \|_{L^4(\Ome; R^{3 \times 3})}
\]
\[
+ \ep^{-1} D_\lambda^{3/4} E_\lambda^{3/4} + \ep V^{1/4} D_\lambda^{3/8} E_\lambda^{5/8}.
\]

\item
 Since $\vc{K}_\lambda = \sigmae \partial_t \vc{u} + \Grad \vc{u} \cdot \left( \vre \vue - \vr \vc{u} \right)$, a direct computation yields
\[
\|\partial_t \vc{K}_\lambda \|_{L^2(\Ome;R^3)} \lesssim \|\partial_t\sigmae \|_{L^2(\Ome)}+ \|\sigmae\|_{L^2(\Ome)}+ \|\partial_t \we \|_{L^2(\Ome;R^3)}+ \|\we \|_{L^2(\Ome;R^3)}.\]
Similarly to the previous step we get
\[
\|\partial_t\sigmae\|_{L^2(\Ome) }\lesssim \|\Grad \we \|_{L^2(\Ome;R^3)} + V^{1/4} \|\Grad \sigmae \|_{L^4(\Ome;R^3) }+\|\sigmae \|_{L^2(\Ome)}+\|\we \|_{L^2(\Ome;R^3)},
\]
and so by \eqref{est1}, \eqref{est2}, \eqref{est3}, and \eqref{est5}
\begin{equation}\label{estI7}
\begin{aligned}
 |I_{7,\lambda}| &\le \ep^4\|\partial_t \vc{K}_\ep\|_{L^2(\Ome;R^3)}\|\partial_t \we \|_{L^2(\Ome;R^3)}
\\ & \lesssim \ep^4 \left( \|\Grad \we \|_{L^2(\Ome; R^{3 \times 3}} + V^{1/4} \|\Grad \sigmae\|_{L^4(\Ome;R^3)}
\right) \|\partial_t \we \|_{L^2(\Ome;R^3)}
\\
& + \ep^4 \left( \|\sigmae\|_{L^2(\Ome)} + \|\we \|_{L^2(\Ome; R^3)} + \|\partial_t \we \|_{L^2(\Ome;R^3)} \right) \|\partial_t \we \|_{L^2(\Ome;R^3)}
\\ &\lesssim E_\lambda + \ep^{2} V^{1/4} \|\Grad \sigmae\|_{L^4(\Ome; R^3)}\sqrt{E_\lambda}.
\end{aligned}
\end{equation}
\end{enumerate}

\medskip

Adding the estimates (\ref{estI1} - \ref{estI7}) we conclude that the inequality \eqref{u12} gives rise to
\begin{equation}\label{u19}
\frac{{\rm d}}{{\rm d}t} E_\lambda +D_\lambda \lesssim
E_\lambda
 +D_{\lambda}^{1/2} E_{\lambda}^{1/2} +
\ep^{-1} D_\lambda^{3/4} E_{\lambda}^{3/4} + \ep V^{1/4} D_{\lambda}^{3/8} E_{\lambda}^{5/8}
\end{equation}
\[
+ V^{\frac{1}{\com 4}} \Big[ {\com V^{\frac{1}{4}}} + \|\Grad \we \|_{L^4(\Ome; R^{3 \times 3})}+\|\Grad \sigmae \|_{L^4(\Ome;R^3)} \Big]
\Big[ D_\lambda^{3/4} E_\lambda^{1/4} + \ep^2 E_\lambda^{1/8} D_\lambda^{3/8} \|\Grad \we \|_{L^4(\Ome; R^{3 \times 3}} \Big]
\]
\[
+ \ep^{2} V^{1/4} \|\Grad \sigmae \|_{L^4(\Ome;R^3)}  \Big[  E_{\lambda}^{1/2} + D_{\lambda}^{1/2} \Big].
\]

\subsubsection{Introducing the modified energy}

Relation (\ref{u19}) contains two integrals that are not controlled by the ``modulated energy'' functional $E_\lambda$. In order to close the
estimates, we differentiate \eqref{u1} in the $x$ variable, multiply the resulting expression by $|\Grad \sigmae |^2 \Grad \sigmae$, and integrate over $\Ome$ obtaining
\begin{equation}\nonumber
\begin{aligned}
\frac{{\rm d}}{{\rm d}t}\int_{\Ome}\frac{1}{4}|\Grad \sigmae |^4 \ \ dx =& -\int_{\Ome} \frac{3}{4}|\Grad \sigmae |^4 \Div \vue +(\Grad \sigmae \cdot \Grad \vue )\cdot(|\Grad \sigmae |^2 \Grad\sigmae )
\\ &+\com|\Grad\sigmae|^2 \Grad\sigmae\cdot \Grad \Div(\vr \we)+(|\Grad\sigmae|^2 \Grad\sigmae)\cdot(\Div\Grad\vue)\sigmae \ \dx.
\end{aligned}
\end{equation}
Applying H\"older's inequality to the terms on the right hand side we get
\[
\frac{{\rm d}}{{\rm d}t}\|\Grad \sigma_\ep\|_{L^4(\Ome;R^3)}^4
\]
\[
\lesssim \|\Grad \sigmae\|_{L^4(\Ome;R^3)}^4 (\|\Grad \we \|_{L^\infty (\Ome; R^{3 \times 3})}  +1 )
\]
\[
+\|\Grad \sigmae \|_{L^4(\Ome;R^3)}^3(
\|\nabla^2_x \we \|_{L^4(\Ome; R^{27}) }+\|\Grad \we \|_{L^4(\Ome; R^{3 \times 3})}+\| \we \|_{L^4(\Ome; R^3)}).
\]
By virtue of \eqref{k3},
\[
\|\Grad \we \|_{L^\infty(\Ome; R^{3 \times 3}) }\lesssim \ep^{-\frac{3}{4}}(\ep \|\nabla^2_x \we \|_{L^4(\Ome; R^{27})}+\|\Grad \we \|_{L^4(\Ome; R^{ 3\times 3})});
\]
whence \eqref{est7} implies
\begin{equation}
\label{u21}
\frac{{\rm d}}{{\rm d}t}\| \Grad \sigmae \|_{L^4(\Ome;R^3)}^2
\end{equation}
\[
\lesssim \|\Grad \sigmae\|_{L^4(\Ome; R^3)}^2 \Big[
\ep^{\frac{1}{4}} \|\nabla^2_x \we \|_{L^4(\Ome; R^{27}) }+ \ep^{-\frac{3}{4}}\|\nabla_x \we \|_{L^4(\Ome; R^{ 3 \times 3}) }+1 \Big]
\]
\[
+ \|\Grad \sigmae\|_{L^4(\Ome; R^3)} \Big[ \|\nabla^2_x \we \|_{L^4(\Ome; R^{27})}+\|\Grad \we \|_{L^4(\Ome; R^{3 \times 3}}+\ep^{-\frac{3}{4}}E_\lambda^{1/2} \Big].
\]

Introducing a new quantity $E_\lambda^* \equiv E_\lambda +\|\Grad \sigmae\|_{L^4(\Ome; R^3)}^2$, and adding
\eqref{u19} to \eqref{u21} we conclude that

\begin{equation}\label{gronwall1}
\frac{{\rm d}}{{\rm d}t} E^*_\lambda +D_\lambda \lesssim
\ep^{-3/4} E^*_\lambda
 +D_{\lambda}^{1/2} (E^*_{\lambda})^{1/2} +
\ep^{-1} D_\lambda^{3/4} (E^*_{\lambda})^{3/4} + \ep V^{1/4} D_{\lambda}^{3/8} (E^*_{\lambda})^{5/8}
\end{equation}
\[
+ V^{\frac{1}{\com 4}} \Big[ {\com V^{\frac{1}{4}}} + \|\Grad \we \|_{L^4(\Ome; R^{3 \times 3})}+ (E^*_\lambda)^{1/2} \Big]
\Big[ D_\lambda^{3/4} (E^*_\lambda)^{1/4} + \ep^2 (E^*_\lambda)^{1/8} D_\lambda^{3/8} \|\Grad \we \|_{L^4(\Ome; R^{3 \times 3})} \Big]
\]
\[
+ \ep^{2} {\com V^{\frac{1}{4}}}(E^*_\lambda)^{1/2}  \Big[  (E^*_{\lambda})^{1/2} + D_{\lambda}^{1/2} \Big]
+  E^*_{\lambda} \Big[
\ep^{\frac{1}{4}} \|\nabla^2_x \we \|_{L^4(\Ome; R^{27}) }+ \ep^{-\frac{3}{4}}\|\nabla_x \we \|_{L^4(\Ome; R^{ 3 \times 3}) } \Big]
\]
\[
+ (E^*_{\lambda})^{1/2} \Big[ \|\nabla^2_x \we \|_{L^4(\Ome; R^{27})}+\|\Grad \we \|_{L^4(\Ome; R^{3 \times 3}}\Big].
\]

\subsubsection{Estimates for $\|\nabla^2_x\vc{w}_\ep\|_{L^4}$ and $ \|\Grad\vc{w}_\ep\|_{L^4} $}

We rewrite \eqref{u3} in the form
\[
- \Div \tn{S} (\Grad \we ) =
-\mu \Delta \we -(\mu/3 + \eta)\Grad \Div \we
\]
\[
=-\vre\partial_t\we -\vre\vue \cdot \Grad  \we -\sigmae \partial_t\vc{u}-(\vre\we +\sigmae \vc{u})\cdot \Grad \vc{u}+\Grad(p(\vr)-p(\vre)),
\]
and denote
\[
G = (p(\vr)-p(\vre)),  \ \vc{F} = -\vre\partial_t\we-\vre\vue \cdot \Grad \we -\sigmae \partial_t\vc{u}-(\vre\we +\sigmae \vc{u})\cdot \Grad \vc{u}.
\]

Then, by the standard elliptic estimates for the Lam\' e system summarized in Section \ref{LM}, we have
\begin{equation}\nonumber
\|\Grad  \we \|_{L^4(\Ome; R^3)}\lesssim \ep^{\frac{1}{4}}\| \vc{F} \|_{L^2(\Ome;R^3) }+\|G\|_{L^4(\Ome)}+\ep^{-1}\|\we \|_{L^4(\Ome;R^3)}
\end{equation}
and
\begin{equation}\nonumber
\|\nabla^2_x \we \|_{L^4(\Ome; R^{27}) }\lesssim \|\vc{F} \|_{L^4(\Ome; R^3)}+ \|\Grad G\|_{L^4(\Ome; R^3) }+\ep^{-2}\|\we \|_{L^4(\Ome; R^3)}.
\end{equation}

Now observe that, by virtue of \eqref{est10},
\begin{equation}\nonumber
\|G\|_{L^4 (\Ome) }\lesssim \|\sigmae \|_{L^4 (\Ome)}
\lesssim d \|\Grad \sigmae\|_{L^4(\Ome; R^3) }+ V^{-\frac{1}{4}}\|\sigmae\|_{L^2(\Ome) }\lesssim \Big( d + V^{-\frac{1}{4}} \Big) (E^*_\lambda)^{1/2},
\end{equation}
and, by {(\ref{est1}-\ref{est5})},
\begin{equation}\nonumber
\| \vc{F} \|_{L^2(\Ome; R^3) }\lesssim \|\partial_t\we \|_{L^2(\Ome; R^3)}+\|\Grad \we \|_{L^2(\Ome; R^{3 \times 3} }+\|\we \|_{L^2(\Ome; R^3) }+\|\sigmae \|_{L^2(\Ome) }\lesssim \ep^{-2} (E^*_\lambda)^{1/2},
\end{equation}
which, {together with \eqref{est7}}, gives rise to
\begin{equation}
\label{u22}
\|\Grad \we \|_{L^4(\Ome; R^{3 \times 3})}\lesssim \Big( \ep^{-\frac{7}{4}} + d + V^{-\frac{1}{4}} \Big) (E^*_\lambda)^{1/2}.
\end{equation}

Finally, in accordance with \eqref{est9}, we have
\[
\|\vc{F} \|_{L^4(\Ome; R^3) }\lesssim \|\partial_t\we \|_{L^4(\Ome;R^3)}+\|\Grad \we \|_{L^4(\Ome;R^3) }+\|\sigmae \|_{L^4(\Ome)}+\|\we \|_{L^4(\Ome; R^3)}
\]
\[
\lesssim {\ep^{-2}}(E^*_\lambda)^{\frac{1}{8}}D_\lambda^{\frac{3}{8}} +  \Big( \ep^{-\frac{7}{4}} + d + V^{-\frac{1}{4}} \Big) (E_\lambda^*)^{1/2},
\]
and, {by \eqref{est10}},
\[
\|\Grad G\|_{L^4 (\Ome; R^3)}\lesssim \|\Grad \sigmae \|_{L^4(\Ome; R^3) }+\|\sigmae\|_{L^4(\Ome)}\lesssim \Big( 1 + d + V^{-1/4} \Big) (E^*_\lambda)^{1/2}.
\]
Combining the above estimates we get
\begin{equation}
\label{u23}
\|\nabla^2_x \we \|_{L^4(\Ome; R^{27}) }\lesssim {\ep^{-2}}(E^*_\lambda)^{\frac{1}{8}}D_\lambda^{\frac{3}{8}} +   \Big( \ep^{-\frac{11}{4}} + d + V^{-\frac{1}{4}} \Big) (E_\lambda^*)^{1/2}.
\end{equation}

Plugging (\ref{u21}), (\ref{u22}) into (\ref{gronwall1}), we obtain

\begin{equation}\label{gronwall}
\frac{{\rm d}}{{\rm d}t} E^*_\lambda +D_\lambda \lesssim
\ep^{-3/4} E^*_\lambda
 +D_{\lambda}^{1/2} (E^*_{\lambda})^{1/2} +
\ep^{-1} D_\lambda^{3/4} (E^*_{\lambda})^{3/4} + \ep V^{1/4} D_{\lambda}^{3/8} (E^*_{\lambda})^{5/8}
\end{equation}
\[
+ V^{\frac{1}{\com 4}} \left[ {\com V^{\frac{1}{4}}} + \left( \ep^{-7/4} + d + V^{-1/4} \right) (E^*_\lambda)^{1/2} \right]
D_\lambda^{3/4} (E^*_\lambda)^{1/4}
\]
\[
+ V^{\frac{1}{\com 4}} \left[ {\com V^{\frac{1}{4}}} + \left( \ep^{-7/4} + d + V^{-1/4} \right) (E^*_\lambda)^{1/2} \right]
\ep^2 (E^*_\lambda)^{1/8} D_\ep^{3/8} \Big( \ep^{-7/4} + d + V^{-1/4} \Big) (E^*_\lambda)^{1/2}
\]
\[
+ \ep^{2} {\com V^{\frac{1}{4}}}(E^*_\lambda)^{1/2}  \Big[  (E^*_{\lambda})^{1/2} + D_{\lambda}^{1/2} \Big]
\]
\[
+ \Big( \ep^{1/4} E^*_\lambda + (E^*_\lambda)^{1/2} \Big) \Big[ {\ep^{-2}}(E^*_\lambda)^{\frac{1}{8}}D_\lambda^{\frac{3}{8}} +  \ep^{-1} \Big( \ep^{-\frac{7}{4}} + d + V^{-\frac{1}{4}} \Big) (E_\lambda^*)^{1/2} \Big].
\]
Finally, applying Young's inequality to the right hand side and regrouping terms of the same order, we may infer that
\begin{equation}
\label{gronwall2}
\frac{{\rm d}}{{\rm d}t}E_\lambda^*+\frac{1}{2}D_\lambda \lesssim
\Big[ E^*_\lambda \left( \ep^{-16/5} + 
V^{-1/4}{\ep^{-1}} \right)
\end{equation}
\[
+  (E^*_\lambda)^3 \left( \ep^{-4} + {\com V \ep^{-7}} \right) + (E^*_\lambda)^{9/5} \left( V^{\com 2/5} \ep^{-12/5} + \ep^{-14/5} + {\com \ep^{16/5} V^{-2/5}} \right)
+ (E^*_\lambda )^{3/2} \left( \ep^{-10/4} + V^{-1/4} \ep^{-3/4} \right) \Big]
\]
as long as the bounds \eqref{h1} hold.

\subsection{Proof of Theorem \ref{TM1}}

With (\ref{gronwall2}) at hand, it is a routine matter to complete the proof of Theorem \ref{TM1}. To begin, we should keep in mind that
(\ref{gronwall2}) holds on condition that (\ref{h1}) is satisfied. Seeing that, in accordance with the embedding (\ref{k3}) and
Poincar\' e inequality (\ref{k4}),
\[
\| \sigmae \|_{L^\infty (\Ome)} + \| \we \|_{L^\infty(\Ome; R^3)}
\]
\[
 \lesssim \ep^{-3/4}
\Big( \| \sigmae \|_{L^4 (\Ome) } + \| \we \|_{L^4 (\Ome; R^3) } +
\ep \| \Grad \sigmae \|_{L^4 (\Ome;R^3)} + \ep \| \Grad \we \|_{L^4 (\Ome; R^{3 \times 3})} \Big)
\]
\[
\lesssim {\com \ep^{3/4}}(\ep + d) \Big( \| \Grad \sigmae \|_{L^4 (\Ome;R^3)} {\com +} \| \Grad \we \|_{L^4 (\Ome; R^{3 \times 3})} \Big) + {\com \ep^{-3/4}}V^{-1/4} (E^*_\lambda)^{1/2}
\]
\[
\lesssim {\com \ep^{-3/4}} V^{-1/4} (E^*_\lambda)^{1/2} + {\com \ep^{-3/4}(\ep+d)} \| \Grad \we \|_{L^4 (\Ome; R^{3 \times 3})},
\]
we may use (\ref{u22}) to deduce that
\[
\| \sigmae \|_{L^\infty (\Ome)} + \| \we \|_{L^\infty(\Ome; R^3)} \lesssim \left( {\com \ep^{-3/4}} V^{-1/4} + \ep^{\com-10/4} \right)  (E^*_\lambda)^{1/2}.
\]
Consequently, for (\ref{h1}) to hold, we need
\bFormula{Z1}
E^*_\lambda (t) \lesssim \min \left\{ \com \ep^{5}, \ep^{3/2} V^{1/2} \right\}, \ t \in [0,T_{\rm max}].
\eF

On the other hand, if (\ref{Z1}) is satisfied, the inequality (\ref{gronwall2}) gives rise to
\bFormula{Z2}
\frac{{\rm d}}{{\rm d}t}E_\lambda^*+\frac{1}{2}D_\lambda \lesssim \Big( {\com \ep^{-16/5}} +  V^{-1/4} \ep^{-1}   \Big) E^*_\lambda.
\eF
Thus, by virtue of Gronwall's lemma,
\bFormula{Z3}
E^*_\lambda (t) \leq E^*_\lambda (0) \exp \left[ C \left( {\com \ep^{-16/5}} +  V^{-1/4} \ep^{-1}  \right) t \right] \ \mbox{for all} \ t \in [0, T_{\rm max}]
\eF
as long as (\ref{Z1}) holds.

Consequently, we may infer that (\ref{Z1}), (\ref{Z3}) holds on the whole time interval $[0,T]$ as soon as
\bFormula{Z4}
E^*_\lambda (0) \leq \exp \left[ - C \left( {\ep^{-16/5}} +  V^{-1/4} \ep^{-1}
 \right) T \right] \min \left\{ \com \ep^{5}, \ep^{3/2}V^{1/2} \right\}.
\eF
In particular, in accordance with Proposition \ref{PM1}, the solution $[\vre, \vue]$ exists on the whole time interval $[0,T]$. Finally, it is easy to check, using in particular (\ref{u3}), that (\ref{MM1}) implies (\ref{Z3}). Theorem \ref{TM1} has been proved.

\section{Applications to problems on thin domains}
\label{A}

Consider the motion of a compressible viscous fluid confined to a thin channel
\[
\Omega_\ep = Q_\ep \times (0,1), \ Q_\ep \subset R^2, \ Q_\ep = \ep Q, \ \ep \to 0,
\]
where $Q \subset R^2$ is a regular planar domain. Furthermore, we assume that the initial density density
$\vr_{0,\ep}$ and the velocity field $\vu_{0,\ep}$ are defined on $\Omega_\ep$, with the integral averages
\[
\frac{1}{|Q_\ep|} \int_{Q_\ep} \vr_{0,\ep}(x_h,y) \ {\rm d}x_h,\ \frac{1}{|Q_\ep|} \int_{Q_\ep} \vr_{0,\ep} \vu_{0,\ep}(x_h,y) \ {\rm d}x_h,\ x_h = (x_1,x_2),
\ y = x_3,
\]
converging weakly (with respect to the $x_3$-variable) to some limit,
\bFormula{Ai1}
\frac{1}{|Q_\ep|} \int_{Q_\ep} \vr_{0,\ep}(x_h,\cdot) \ {\rm d}x_h \to \vr_0 ,\
\frac{1}{|Q_\ep|} \int_{Q_\ep} \vr_{0, \ep} \vu_{0,\ep}(x_h,\cdot) \ {\rm d}x_h \to (\vr \vu)_0 \ \mbox{weakly in} \ L^1(0,1)
\eF
as $\ep \to 0$.

We suppose that $\vr_\ep = \vr_\ep(t,x)$, $\vu_\ep = \vu_\ep(t,x)$ is a solution of the compressible Navier-Stokes system (\ref{i1} - \ref{i3}),
supplemented with the slip conditions (\ref{i4}), and the initial data $[\vr_{0,\ep}, \vu_{0,\ep}]$. As the limit data depend only on the
$x_3 \equiv y$-variable, a candidate for the limit problem is the $1D$ compressible Navier-Stokes system:
\bFormula{Ai4}
\partial_t \vr + \partial_y (\vr u) = 0,
\eF
\bFormula{Ai5}
\partial_t (\vr u) + \partial_y (\vr u^2) + \partial_y p(\vr) = \nu \partial^2_{y,y} u ,\ \nu = \frac{4}{3} \mu + \eta,
\eF
where $\vr = \vr(t,y)$, $u= u(t,y)$. Indeed, we claim the following result, see \cite[Theorem 2.1]{BeFeNov}:

\bTheorem{A1}
Let $Q \subset R^2$ be a bounded Lipschitz domain. Suppose that the pressure $p = p(\vr)$ satisfies
\[
p \in C[0, \infty) \cap C^2(0, \infty), \ p'(\vr) > 0 \ \mbox{for all}\ \vr > 0,\
\lim_{\vr \to \infty} \frac{p'(\vr)}{\vr^{\gamma - 1}} = p_\infty > 0, \ \gamma > \frac{3}{2},
\]
and that
the viscous stress tensor $\tn{S}$ is given by (\ref{i3}), with the viscosity coefficients
\[
\mu > 0, \ \eta > 0.
\]

Let
\[
\frac{1}{|Q_\ep|} \int_{Q_\ep} \vr_{0,\ep}(x_h,\cdot) \ {\rm d}x_h \to \vr_0 ,\
\frac{1}{|Q_\ep|} \int_{Q_\ep} \vr_{0, \ep} \vu_{0,\ep}(x_h,\cdot) \ {\rm d}x_h \to (\vr \vc{u})_0 \ \mbox{weakly in} \ L^1(0,1),
\]
\[
\frac{1}{|Q_\ep|} \intOe{ \left[ \frac{1}{2} \vr_{0,\ep} |\vu_{0,\ep}|^2 + H(\vr_{0,\ep}) \right] } \to
\int_0^1  \left[ \frac{1}{2 \vr_{0}} |(\vr u)^3_0 |^2 + H(\vr_{0}) \right] \ {\rm d}y
\]
where $\vr_0 > \underline{\vr}$, $u_0 \equiv (\vr u)^3_0 / \vr_0$ satisfy
\bFormula{A10a}
\left\{
\begin{array}{c}
\vr_0 \in C^{1 + \beta}[0,1], \ u_0  \in C^{2 + \beta}[0,1],\ \beta > 0,\\ \\
\mbox{with the compatibility conditions}\  u_0|_{y=0,1} = \partial^2_{y,y} u_0|_{y = 0,1} = \partial_y \vr_0|_{0,1} = 0.
\end{array} \right\}
\eF
Let $[\vr_\ep, \vu_\ep]$ be a \emph{finite energy weak} solution of the barotropic Navier-Stokes system (\ref{i1} - {\com \ref{i4}})
in $(0,T) \times \Omega_\ep$, emanating from the initial data
\[
\vr_\ep (0,\cdot) = \vr_{0,\ep},\ (\vr_\ep \vu_\ep )(0, \cdot) = \vr_{0,\ep} \vu_{0, \ep}.
\]

Then
\[
{\rm ess} \sup_{t \in (0,T)} \frac{1}{|Q_\ep|} \int_0^1 \int_{Q_\ep}  \left| \vr_\ep  (t,x_h,y) - \vr (t, y) \right|^\gamma
\ {\rm d}x_h \ {\rm d}y \to 0,
\]
and
\[
{\rm ess} \sup_{t \in (0,T)} \frac{1}{|Q_\ep|} \int_0^1 \int_{Q_\ep} \left| {\vr_\ep}{\vu_\ep}(t,x_h,y) - [0,0,\vr u](t, y) \right|^{2\gamma/(\gamma + 1)}
\ {\rm d}x_h \ {\rm d}y \to 0
\]
as $\ep \to 0$,
where $[\vr, u]$ is the unique solution of the $1D$ Navier-Stokes system (\ref{Ai4}), (\ref{Ai5}), with the initial data $[\vr_0, u_0]$
satisfying the no-slip conditions
\[
v(t,0) = v(t,1) = 0.
\]
\eT

Adapting Theorem \ref{TM1} to the present situation we obtain the following result:

\Cbox{Cgrey}{

\bTheorem{Z1}
In addition to the hypotheses of Theorem \ref{TA1}, suppose that $Q$ is a uniformly $C^4-$domain and that the initial data
$[\vr_{0,\ep}, \vu_{0, \ep}]$ belong to the regularity class specified in Theorem \ref{TM1} and satisfy
\[
\left\| \vr_{0,\ep} - \vr_0 \right\|_{W^{1,4}(\Omega_\ep)} + \left\| \vu_{0,\ep} - [0,0,u_0] \right\|_{W^{2,2}(\Omega_\ep;R^3)} \leq
\ep^{\com 5} \exp \left( - C \ep^{\com -16/5}  T \right).
\]

Then $[\vr_\ep, \vu_\ep]$ is a strong solution in the class specified in Proposition \ref{PM1}.

\eT

}

\bRemark{Z1}

Note that, in accordance with the general weak-strong uniqueness property established in \cite{FeNoJi}, the weak and strong solution emanating from the same initial data coincide as long as the latter exists.

\eR

\bRemark{Z2}

Since the initial data for the limit problem satisfy (\ref{A10a}), the solution $[\vr, u]$ is smooth as required in Theorem \ref{TM1},
see Kazhikhov \cite{KHA2}.

\eR

\def\cprime{$'$} \def\ocirc#1{\ifmmode\setbox0=\hbox{$#1$}\dimen0=\ht0
  \advance\dimen0 by1pt\rlap{\hbox to\wd0{\hss\raise\dimen0
  \hbox{\hskip.2em$\scriptscriptstyle\circ$}\hss}}#1\else {\accent"17 #1}\fi}

\end{document}